\documentclass[12pt,3p,nopreprintline]{elsarticle}
\pdfoutput=1
\ifpdf\fi

\usepackage[T1]{fontenc}
\usepackage[slantedGreek]{newtx}
\usepackage{microtype}

\usepackage{amsmath}

\usepackage[]{hyperref}
\newcommand{\dd}{\partial}

\begin{document}

\begin{frontmatter}

\title{Derivatives of the full QR factorisation and of \\ the factored-form and compact WY representations}
\author{Stefanos-Aldo Papanicolopulos}
\ead{S.Papanicolopulos@ed.ac.uk}
\address{School of Engineering, Institute for Infrastructure and Environment, The University of Edinburgh,\\ The King's Buildings, Edinburgh, United Kingdom}

\begin{abstract}
QR factorisation plays an important role in matrix computations. Within the context of optimisation and of automatic differentiation of such computations, we need to compute the derivative of this factorisation. For tall matrices, however, existing results only cover the so-called thin case.

We provide for the first time expressions for the derivative of the full QR factorisation of a tall matrix, in the usual case where the $Q$ factor is a product of Householder reflections. These expressions are obtained based on novel results for the derivative of the compact WY representation of $Q$, which also yield the derivative of the factored-form representation of $Q$, both of which are useful on their own.
These three results can be used directly in applications such as variable projection for solving separable non-linear least squares problems, and can also extend the current linear algebra capabilities of automatic differentiation frameworks.

\end{abstract}

\begin{keyword}
QR factorisation, QR decomposition, compact WY representation, factored-form representation, differentiation, automatic differentiation
\MSC[2010] 65F25 \sep             65D25       \end{keyword}
\end{frontmatter}

\section{Introduction}

Matrix computations form a cornerstone of linear algebra~\citep[see e.g.][]{GolubVanLoan2013}. Among these computations, matrix factorisations (or decompositions) such as Cholesky, LU, QR and SVD (singular value decomposition) play an important role and have therefore been extensively studied both theoretically and from an implementation point of view. Robust and efficient implementations are provided in software libraries such as LAPACK~\cite{LAPACK}.

Expressions for the derivatives of matrix computations, including matrix factorisations, are needed when evaluating gradients of quantities involving such computations, either within the framework of automatic differentiation or more generally in symbolic calculations. Useful results have been published for the differentiation of matrix computations~\cite{Minka2000, Giles2008inbook}, and more specifically for matrix factorisations (Cholesky~\cite{Murray2016,Seeger2017}, LU~\cite{deHoog2011}, QR~\cite{Walter2012,Seeger2017,Townsend2018,roberts2020qr}, and SVD~\cite{Townsend2016, Seeger2017, Kanchi2025}).

We focus here on the QR factorisation (also known as QR decomposition), which plays an important role for example in solving the linear least squares problem~\citep{GolubVanLoan2013}. For tall matrices, there are known results for the derivative of the so-called thin QR factorisation~\citep{Walter2012,roberts2020qr}, with early related results dating back to the development of variable projection~\cite{GolubPereyra1973}. This derivative is implemented for example in the gradient calculations in frameworks such as PyTorch and JAX. These frameworks however do not implement the full factorisation,\footnote{This has been requested for different frameworks, see \url{https://github.com/jax-ml/jax/issues/22264}, \url{https://github.com/tensorflow/tensorflow/issues/6504\#issuecomment-832180398}, \url{https://discuss.pytorch.org/t/200312}}
and indeed we are not aware of complete theoretical results for this case. A reason given for this is that while the thin case does not depend on the QR algorithm used, the full case does\footnote{\url{https://github.com/google/jax/issues/10450\#issuecomment-1109982015}}, a fact already appreciated in early work by~\citet{Kaufman1975}.

To address this dependence on the algorithm used, in this paper we consider the compact WY representation for the $Q$ factor~\citep{Schreiber1989}, present formulas for its derivative, and use them to obtain expressions for the derivative of the full QR factorisation and also for the derivative of the factored-form representation of $Q$.
These are novel results, that fill a gap in the relevant literature and can be used to complete the implementation of automatic differentiation of the QR factorisation, as well as for symbolic computations.

Section~\ref{sec:notation} introduces the notation used, while sec.~\ref{sec:fullthinQR} provides a brief overview of the QR factorisation. Section~\ref{sec:thinderiv} summarises known results for the derivative of the QR factorisation, and provides further information on previous attempts to compute derivatives of the full factorisation. New results are provided in sec.~\ref{sec:cWYderiv} for the derivative of the compact WY representation and of the factored-form representation, and in sec.~\ref{sec:fullderiv} for the derivative of the full QR factorisation. We discuss further the importance of these results in section~\ref{sec:contribution}.

\section{Notation}
\label{sec:notation}

We consider the following partitions of a square $m \times m$ matrix $A_{mm}$
\begin{equation}
    A_{mm} = \begin{bmatrix} A_{mn} & A_{mp} \end{bmatrix}
           = \begin{bmatrix} A_{nm} \\ A_{pm} \end{bmatrix}
           = \begin{bmatrix} A_{nn} & A_{np} \\ A_{pn} & A_{pp} \end{bmatrix}
\end{equation}
where $p=m-n \geq 0$. The indices indicate the dimensions of each matrix, therefore $A_{mn}$ is a tall matrix, while $A_{nm}$ is a wide one, except for the case $p=0$ where both are square and equivalent to $A_{nn}$. We consider here as ``tall'' any matrix with more rows than columns, without any requirements on the ratio of the two being large.

In this notation, the name of the indices is important: for example, $A_{pp}$ is not simply a $p \times p$ matrix, but it is specifically the bottom-right $p \times p$ block of the $m \times m$ matrix $A_{mm}$. On the other hand, we can consider for example the $A_{mn}$ matrix (and its partition into $A_{nn}$ and $A_{pn}$ blocks) without the need to introduce explicitly an $A_{mm}$ matrix and its $A_{mp}$ block.

To work with upper triangular matrices we introduce the $U$ matrix, whose elements are one above or on the diagonal and zero elsewhere. Similarly, to work with strictly lower triangular matrices, we introduce the $\hat{L}$ matrix, whose elements are one below the diagonal and zero elsewhere.
For example, in the $4 \times 4$ case we have
\begin{equation}
  U = \begin{bmatrix}
      1 & 1 & 1 & 1 \\
      0 & 1 & 1 & 1 \\
      0 & 0 & 1 & 1 \\
      0 & 0 & 0 & 1
  \end{bmatrix}    
,\qquad
  \hat{L} = \begin{bmatrix}
      0 & 0 & 0 & 0 \\
      1 & 0 & 0 & 0 \\
      1 & 1 & 0 & 0 \\
      1 & 1 & 1 & 0
  \end{bmatrix}    
\; .
\end{equation}
The $U$ and $\hat{L}$ matrices are only involved in element-wise products, such as $U \circ A_{mn}$, so their dimensions are always implied by the second term of the product.

We will be considering derivatives of matrices with respect to a given scalar variable, such as $\partial A_{mn} / \partial \alpha$. The actual variable $\alpha$ is not important in this paper, so for simplicity we will write the derivative $\partial A_{mn} / \partial \alpha$ as $\dd A_{mn}$. The results in this paper can therefore also be considered from the viewpoint of matrix differentials, or of matrix perturbations.

\section{The QR factorisation}
\label{sec:fullthinQR}

We consider in this paper only real matrices, although it should be easy to extend the results to the complex case. A matrix $A_{mn}$ can be written as the (\emph{full}) QR factorisation \citep[sec~5.2]{GolubVanLoan2013}
\begin{equation}\label{eq:QR}
    A_{mn} = Q_{mm} R_{mn}
\end{equation}
where $Q_{mm}$ is orthogonal and $R_{mn}$ is upper triangular. Partitioning~\eqref{eq:QR}, we obtain the \emph{thin} (or \emph{reduced}) QR factorisation
\begin{equation}\label{eq:thinQR}
    A_{mn} = \begin{bmatrix} Q_{mn} & Q_{mp} \end{bmatrix} \begin{bmatrix} R_{nn} \\ 0 \end{bmatrix}
           = Q_{mn}  R_{nn}
\end{equation}

The QR factorisation is generally not unique. If $A_{mn}$ has full column rank, there is a unique thin QR factorisation for which the diagonal elements of $R_{nn}$ are all positive. The $Q_{mp}$ block of the full factorisation has orthonormal columns and is orthogonal to $Q_{mn}$, but is not uniquely defined.

A widely used class of algorithms for computing the QR factorisation expresses $Q_{mm}$ as the product of $n$ Householder reflections, i.e.\ symmetric orthogonal matrices of the form
\begin{equation}\label{eq:householder}
    H_{mm}^{(i)} = I - \tau^{(i)} v_m^{(i)} \bigl(v_m^{(i)}\bigr)^T  \qquad i=1 \ldots n
\end{equation}
where each $\tau^{(i)}$ is a scalar and each $v_m^{(i)}$ is a vector with $m$ elements, of which the first $i-1$ are zero and the $i$-th is one.
Householder QR algorithms generally return the \emph{factored-form} representation of $Q_{mm}$, which consists of a lower unit triangular matrix $Y_{mn}$ whose $i$-th column is the vector $v_m^{(i)}$, and of a vector containing the $\tau^{(i)}$ coefficients. The $\tau^{(i)}$ coefficients are stored for efficiency, though they can be computed from the $v_m^{(i)}$ vectors, based on the orthogonality of $H_{mm}^{(i)}$, as
\begin{equation}
    \tau^{(i)} = \frac{2}{\bigl(v_m^{(i)}\bigr)^T v_m^{(i)} }
\end{equation}

\Citet{Schreiber1989} expressed the product of Householder reflections giving the $Q_{mm}$ factor using the \emph{compact WY} representation
\begin{equation}\label{eq:compactWY}
    Q_{mm} = I_{mm} - Y_{mn} T_{nn} Y_{mn}^T
\end{equation}
where $Y_{mn}$ is the lower unit triangular matrix also appearing in the factored-form representation, and $T_{nn}$ is upper triangular with its diagonal elements being the $\tau^{(i)}$ coefficients. This representation can result in improved computational efficiency, e.g.\ in the recursive algorithm of \citet{Elmroth2000}.
The $T_{nn}$ matrix can be (inneficiently) computed directly from $Y_{mn}$, since the orthogonality of $Q_{mm}$ yields
\begin{equation}\label{eq:Tnn_inv}
    T_{nn}^{-1} + \bigl( T_{nn}^{-1} \bigl)^T = Y_{mn}^T Y_{mn}
    \;.
\end{equation}

An advantage of the factored-form and compact WY representations is that they allow for more efficient computation of matrix products involving the $Q$ factor (or one of its blocks), through the use of dedicated matrix multiplication routines (ORMQR and GEMQRT, respectively, in LAPACK).

Expressing $Q_{mm}$ through the use of the compact WY representation (and, more generally, as the product of Householder reflections) introduces a unique definition of the $Q_{mp}$ block, which we will use in the following to obtain the derivative of the entire $Q_{mm}$ factor. There are however algorithms (e.g.\ using Givens rotations) that will yield $Q_{mm}$ factors that cannot be represented using the compact WY representation, as shown in a simple example in~\ref{sec:minimalexample}.

Routines that calculate Householder reflections, such as the LARFG and LARFGP routines in LAPACK or the \texttt{house} routine in~\citep[Algorithm~5.1.1]{GolubVanLoan2013}, generally allow for a special case where $\tau^{(i)}=0$ for some $i$. This however introduces a discontinuity, therefore in the following we will assume that $\tau^{(i)} \neq 0$ for all reflections, which means that $T_{nn}$ is invertible.

\section{Known results on the derivative of the QR factorisation}
\label{sec:thinderiv}
\label{sec:thinQRderivative}

Given the derivative $\dd A_{mn}$, we want to calculate the derivatives $\dd Q_{mm}$ and $\dd R_{nn}$ of its full QR factorisation (obviously $\dd R_{pn}=0)$. We assume that all required derivatives do exist; this is not necessarily the case for a given algorithm, see e.g.\ \cite{Coleman1984} for a discussion on the continuity of $Q_{mp}$.

\subsection{Derivative of the thin QR factorisation}

We present in this section known results for the thin factorisation~\citep[see e.g.][]{Walter2012,roberts2020qr}. 
We assume that $A_{mn}$ has full column rank, in which case $R_{nn}$ is invertible. 
Differentiating the thin factorisation~\eqref{eq:thinQR} and multiplying on the right with $R_{nn}^{-1}$, we get
\begin{equation}\label{eq:Bmn}
    \underbrace{(\dd A_{mn}) R_{nn}^{-1}}_{B_{mn}}
    =
    \dd Q_{mn}
    +
    Q_{mn}
    \underbrace{(\dd R_{nn}) R_{nn}^{-1}}_{\Psi_{nn}}
\end{equation}
and, multiplying on the left with $Q_{mn}^T$,
\begin{equation}\label{eq:Enn}
    \underbrace{Q_{mn}^T B_{mn}}_{E_{nn}}
    =
    \underbrace{Q_{mn}^T \dd Q_{mn}}_{\Omega_{nn}} {} + \Psi_{nn}
    \;,
\end{equation}
where we use underbraces to introduce new quantities.

We see that $\Omega_{nn}$ is skew symmetric (by differentiating the orthonormality condition $Q_{mn}^T Q_{mn} = I_{nn}$) and $\Psi_{nn}$ is upper triangular; therefore
\begin{align}
    \Psi_{nn} &= (U \circ E_{nn}) + (\hat{L} \circ E_{nn})^T \label{eq:Psinn}
    \;.
\end{align}
Equation~\eqref{eq:Bmn} then yields the derivative of the thin QR factorisation
\begin{subequations}\label{eq:thinderiv}
\begin{align}
    \dd R_{nn} &= \Psi_{nn} R_{nn} \label{eq:dRnn} \\
    \dd Q_{mn} &= B_{mn} - Q_{mn} \Psi_{nn} \label{eq:dQmn}
    \;.
\end{align}
\end{subequations}

This way of calculating the derivative of the thin QR factorisation only involves the thin factors, resulting in more efficient expressions for calculating the derivative.

\subsection{Alternative form using skew symmetric matrix}

Differentiating the orthogonality condition $Q_{mm}^T Q_{mm} = I_{mm}$ it is easy to see that
\begin{equation}\label{eq:dQmm}
    \dd Q_{mm} = Q_{mm} \Omega_{mm}
\end{equation}
where $\Omega_{mm}$ is a skew symmetric matrix. The (also skew symmetric) $\Omega_{nn}$ block of $\Omega_{mm}$ has already appeared in~\eqref{eq:Enn}, from which it can be calculated as
\begin{equation}
    \label{eq:Omegann}
    \Omega_{nn} = (\hat{L} \circ E_{nn}) - (\hat{L} \circ E_{nn})^T
\end{equation}
Multiplying~\eqref{eq:Bmn} on the left with $Q_{mp}^T$ yields
\begin{equation}
    \label{eq:Omegapn}
        \Omega_{pn} = Q_{mp}^T B_{mn}
\end{equation}

Equation~\eqref{eq:dQmm} therefore gives the alternative form
\begin{equation}\label{eq:dQmn_Omega}
    \dd Q_{mn} = Q_{mm} \Omega_{mn}
\end{equation}
where $\Omega_{mn}$ is given by~\eqref{eq:Omegann} and~\eqref{eq:Omegapn}. This form provides a more structured form for the derivative $\dd Q_{mn}$, which can be useful in symbolic calculations, but is less efficient to calculate than~\eqref{eq:dQmn}.

\subsection{Known results for the full factorisation}

Equation~\eqref{eq:dQmm} also gives
\begin{equation}\label{eq:dQmp_Omega}
    \dd Q_{mp} = Q_{mm} \Omega_{mp}
\end{equation}
where $\Omega_{np} = -\Omega_{pn}^T$ and is therefore known from~\eqref{eq:Omegapn}, but the skew symmetric block $\Omega_{pp}$ depends on the way the QR factorisation is calculated.
\citet{Kaufman1975} presented this result and made the conjecture that there exists a $Q_{mp}$ matrix for which $\Omega_{pp} = 0$, but we are not aware of any proof of this conjecture. For the $Q_{mp}$ factor obtained using Householder reflections, it is certainly the case that $\Omega_{pp} \neq 0$, as we prove in sec.~\ref{sec:fullderiv}.

More recently, algorithm~1 in \citep{Walter2012} similarly assumes that $\Omega_{pp}=0$ and therefore does not give the correct derivative of the full QR factorisation. Equation~(A16) in~\cite{Baerligea2023} also gives an explicit expression for $\dd Q_{mp}$, independent of the algorithm used for obtaining $Q_{mp}$, which is however incorrect.

We see therefore that the calculation of $\dd Q_{mp}$, and therefore of the full QR factorisation, is still an open question. To address this, knowing that $Q_{mp}$ depends on the algorithm used, we consider next the derivative of the compact WY representation of the $Q_{mm}$ factor.

\section{Derivative of the compact WY representation}
\label{sec:cWYderiv}

The $Q_{nn}$ block of the compact WY representation~\eqref{eq:compactWY} is
\begin{equation}\label{eq:compactWYnn}
    Q_{nn}
    = I_{nn} - Y_{nn} T_{nn} Y_{nn}^T
    = I_{nn} + Y_{nn} ( \underbrace{-T_{nn} Y_{nn}^T}_{S_{nn}} )
\end{equation}
therefore its derivative is
\begin{equation}\label{eq:dQnn}
    \dd Q_{nn} = ( \dd Y_{nn} ) S_{nn} + Y_{nn} \dd S_{nn}
    \;.
\end{equation}
Since $Y_{nn}$ is unit lower triangular, it is always invertible. We have also assumed that $T_{nn}$ is invertible, therefore $S_{nn}$ is also invertible. Equation~\eqref{eq:dQnn} then can be written as
\begin{equation}\label{eq:Cnn}
    \underbrace{ Y_{nn}^{-1} (\dd Q_{nn}) S_{nn}^{-1} }_{C_{nn}}
    =
    Y_{nn}^{-1} \dd Y_{nn}
    + \dd S_{nn} S_{nn}^{-1}
    \;.
\end{equation}
The r.h.s.\ of~\eqref{eq:Cnn} is the sum of a strictly lower triangular term (as the diagonal of $\dd Y_{nn}$ is zero) and an upper triangular term, therefore
\begin{align}
    \dd Y_{nn} &= Y_{nn} (\hat{L} \circ C_{nn})  \label{eq:dYnn} \\
    \dd S_{nn} &= (U \circ C_{nn}) S_{nn}        \label{eq:dSnn}
    \;.
\end{align}

Differentiating the $Q_{pn}$ block of the compact WY representation~\eqref{eq:compactWY} yields
\begin{equation}    \dd Q_{pn} = ( \dd Y_{pn} ) S_{nn} + Y_{pn} \dd S_{nn}
\end{equation}
therefore
\begin{equation}\label{eq:dYpn}
    \dd Y_{pn} = (\dd Q_{pn}) S_{nn}^{-1}  -  Y_{pn} (U \circ C_{nn})
    \;.
\end{equation}
Equations~\eqref{eq:dYnn} and~\eqref{eq:dSnn}, together with the definition of $S_{nn}$ in~\eqref{eq:compactWYnn}, yield the derivative of $T_{nn}$ as
\begin{equation}\label{eq:dTnn}
    \dd T_{nn} =
    (U \circ C_{nn}) T_{nn}
    - T_{nn} (\hat{L} \circ C_{nn} )^T 
    \;.
\end{equation}

Equations~\eqref{eq:dYnn}, \eqref{eq:dYpn} and~\eqref{eq:dTnn} provide the derivative of the compact WY representation, but depend on having already calculated the derivative $\dd Q_{mn}$, using for example equation~\eqref{eq:dQmn} or~\eqref{eq:dQmn_Omega}.
Substituting and simplifying we obtain
\begin{equation}
    C_{nn}         = \underbrace{Y_{nn}^{-1} (B_{nn} - \Psi_{nn}) S_{nn}^{-1}}_{C^*_{nn}} - S_{nn} \Psi_{nn} S_{nn}^{-1}
\end{equation}
and the derivative of the compact WY representation can then be written as
\begin{subequations}
\label{eq:diffCWY}
\begin{align}
    \dd Y_{nn} &= Y_{nn} (\hat{L} \circ C^*_{nn})  \label{eq:dYnn2} \\
    \dd Y_{pn} &= B_{pn} S_{nn}^{-1} - Y_{pn} (U \circ C^*_{nn}) \\
    \dd T_{nn} &= (U \circ C^*_{nn}) T_{nn} - T_{nn} (\hat{L} \circ C^*_{nn} )^T + S_{nn} \Psi_{nn} Y_{nn}^{-T}
    \;.
\end{align}
\end{subequations}

Equations~\eqref{eq:diffCWY} also provide the derivative of the factored-form representation of $Q_{mm}$, since the $Y_{mn}$ matrix is the same and $\tau^{(i)}$ are the diagonal elements of $T_{nn}$. Indeed, it is easy to see that
\begin{equation}
    \dd \tau^{(i)} = \bigl( C^*_{nn}(i,i) - \Psi_{nn}(i,i) \bigr) \tau^{(i)}
\end{equation}
Since the factored-form representation does not include the $T_{nn}$ matrix, the $T_{nn}^{-1}$ matrix needed to compute $S_{nn}^{-1}$  is obtained in this case using eq.~\eqref{eq:Tnn_inv}.

\section{Derivative of the full QR factorisation}
\label{sec:fullderiv}

Since expressions for $\dd R_{mn}$ and $\dd Q_{mn}$ are already known, calculating the derivative of the full QR factorisation only requires the additional calculation of $\dd Q_{mp}$. Differentiating the $Q_{mp}$ block of the compact WY representation~\eqref{eq:compactWY}, we get
\begin{equation}\label{eq:dQmp_init}
    \dd Q_{mp}
    = \dd ( I_{mp} - Y_{mn} T_{nn} Y_{pn}^T )
    = - (\dd Y_{mn}) T_{nn} Y_{pn}^T - Y_{mn} (\dd T_{nn}) Y_{pn}^T - Y_{mn} T_{nn} (\dd Y_{pn}^T)
\end{equation}
which can be computed, since the necessary derivatives were obtained in section~\ref{sec:cWYderiv}.
Indeed, using eqs.~\eqref{eq:diffCWY}, we get after some calculations 
\begin{equation}\label{eq:dQmp_YS}
    \dd Q_{mp} =
    \dd Q_{mn} ( Y_{pn} Y_{nn}^{-1} )^T
       - Y_{mn} T_{nn} S_{nn}^{-T} (\dd Q_{pn} - \underbrace{Y_{pn} Y_{nn}^{-1}}_{Z_{pn}} \dd Q_{nn})^T 
    \;.
\end{equation}

Using the orthonormality condition $Q_{mn}^T Q_{mn} = I_{nn}$ we can show that
\begin{equation}
    Y_{mn} T_{nn} S_{nn}^{-T} = Q_{mn} + Q_{mp} Z_{pn}
\end{equation}
therefore eq.~\eqref{eq:dQmp_YS} becomes
\begin{equation}\label{eq:dQmp}
    \dd Q_{mp} =
       (\dd Q_{mn}) Z_{pn}^T
     - (Q_{mn} + Q_{mp} Z_{pn}) ( \dd Q_{pn} - Z_{pn} \dd Q_{nn} )^T
    \;.
\end{equation}
It is also easy to show that
\begin{equation}\label{eq:Zpn}
    Z_{pn} = Q_{pn} (Q_{nn} - I_{nn})^{-1}
\end{equation}
therefore it is possible to write $\dd Q_{mp}$ in~\eqref{eq:dQmp} purely as a function of $\dd Q_{mn}$ and $Q_{mm}$, without direct reference to the form in which $Q_{mm}$ is computed (in this case using Householder reflections). This does not mean however that eqs.~\eqref{eq:dQmp} and~\eqref{eq:Zpn} yield $\dd Q_{mp}$ independently of the form of $Q_{mm}$, as in general they will not be valid when $Q_{mm}$ is not the product of Householder reflections. This is demonstrated in~\ref{sec:minimalexample} through a simple example.

From eq.~\eqref{eq:dQmp} we can also obtain the $\Omega_{pp}$ block needed in eq.~\eqref{eq:dQmp_Omega} as
\begin{equation}\label{eq:Omegamp}
    \Omega_{pp} = \Omega_{pn} Z_{pn}^T - Z_{pn} \Omega_{pn}^T - Z_{pn} \Omega_{nn} Z_{pn}^T
    \;.
\end{equation}
which shows that setting $\Omega_{pp}=0$ does not give the correct derivative $\dd Q_{mp}$ when using Householder reflections.

\section{Contribution}
\label{sec:contribution}

The main contribution of this paper are three novel results on the derivation of the $Q$ factor of the QR factorisation, by providing explicit expressions for the derivatives of a) the compact WY representation; b) the factored-form representation; and c) the $Q_{mp}$ block (when $Q_{mm}$ is the product of Householder reflections). These address a significant gap, as we are not aware of similar results in the literature.

These results can be explicitly used in computations. As an example, $\dd Q_{mp}$ appears in the approach of \citet{Kaufman1975} for solving separable non-linear least squares problems using variable projection, where an approximation is used instead. Using the exact derivative should accelerate convergence, and possibly avoid stalling of the algorithm~\cite{OLeary2013}, while still requiring fewer matrix multiplications compared to the original approach of~\citet{GolubPereyra1973}.

Our results are also important within the context of automatic differentiation by
providing the derivative of the full factor $Q_{mm}$, a missing capability that (as mentioned in the introduction) has been requested for different automatic differentiation frameworks. These results also provide the derivatives of factorisation routines returning the factored-form representation (e.g.\ GEQRF) or the compact WY representation (e.g.\ GEQRT3). Such functions are often used, together with appropriate matrix multiplication routines, to reduce computational cost compared to explicit computation of the $Q$ factor.

In all the above cases, the contribution consists in introducing new, previously unavailable capabilities; therefore there is no scope for comparing computational performance with existing implementations. While we have not specifically addressed the computational cost of the proposed expressions, we have sought to provide simple forms, that should allow for efficient implementations based on optimised linear algebra libraries.

\section*{Acknowledgements}

This research was supported through a sabbatical leave from the School of Engineering of the University of Edinburgh.

\setlength{\bibsep}{0pt plus 1pt}
\bibliographystyle{elsarticle-num-names}
\bibliography{qr}

\appendix

\section{Minimal example of calculations}\label{sec:minimalexample}

The smallest tall matrix is $2 \times 1$. This provides a minimal, but still meaningful, example of the calculations presented in this paper.
Consider the generic (real) matrix
\begin{equation}
    A_{mn} = \begin{bmatrix} a_{11} \\ a_{21} \end{bmatrix}
\end{equation}
In this case $R_{nn}$ has a single, diagonal element. Requiring this element to be positive easily gives the unique thin factorisation
\begin{equation}\label{eq:QRthin21}
    Q_{mn} = \frac{1}{r_{11}} \begin{bmatrix} a_{11} \\ a_{21} \end{bmatrix}
    ,\quad
    R_{nn} = \begin{bmatrix} r_{11} \end{bmatrix}
    ,\quad \text{with }
    r_{11} = \sqrt{a_{11}^2 + a_{21}^2}
    \;.
\end{equation}

Given the derivative
\begin{equation}
    \dd A_{mn} = \begin{bmatrix} \dd a_{11} \\ \dd a_{21} \end{bmatrix}
\end{equation}
we can use~\eqref{eq:QRthin21} to calculate the QR factors $Q_{mn}^*$ and $R_{nn}^*$  of $A_{mn}^* = A_{mn} + \epsilon \dd A_{mn}$. We can then calculate directly the derivatives $\dd Q_{mn}$ and $\dd R_{nn}$. For the latter we get
\begin{equation}
    \dd r_{11} = \lim_{\epsilon \to 0} \frac{r_{11}^* - r_{11}}{\epsilon}
    = \cdots = \frac{ a_{11} \dd a_{11} + a_{21} \dd a_{21} }{ r_{11} }
\end{equation}
while for the former we similarly get
\begin{equation}
    \dd Q_{mn} = \frac{ a_{11} \dd a_{21} - a_{21} \dd a_{11} }{r_{11}^3}
                 \begin{bmatrix} -a_{21} \\ a_{11} \end{bmatrix}
\end{equation}
It is easy to check that eqs.~\eqref{eq:thinderiv} give the same result.

For the $Q_{mp}$ block, there are in this case two possible cases
\begin{equation}
    Q_{mp}^h = \frac{1}{r_{11}} \begin{bmatrix} a_{21} \\ -a_{11} \end{bmatrix}
    ,\quad
    Q_{mp}^g = -Q_{mp}^h
\end{equation}
with derivatives directly calculated (using limits) as
\begin{equation}
    \dd Q_{mp}^h = \frac{ a_{11} \dd a_{21} - a_{21} \dd a_{11} }{r_{11}^3}
                    \begin{bmatrix} a_{11} \\ a_{21} \end{bmatrix}
    ,\quad
    \dd Q_{mp}^g = - Q_{mp}^h
\end{equation}

Using $Q_{mp}^h$ results in a $Q_{mm}$ factor that is a Householder reflection, which can be directly expressed in a compact WY representation
\begin{equation}\label{eq:Qmmh}
    Q_{mm}^h = \frac{1}{r_{11}} \begin{bmatrix} a_{11} & a_{21} \\ a_{21} & -a_{11} \end{bmatrix}
    = \begin{bmatrix} 1 & 0 \\ 0 & 1 \end{bmatrix}
      - \begin{bmatrix} 1 \\ y_{21} \end{bmatrix}
      \begin{bmatrix} t_{11} \end{bmatrix}
      \begin{bmatrix} 1 & y_{21} \end{bmatrix}
\end{equation}
with
\begin{equation}\label{eq:Qmmh_yt}
    y_{21} = \frac{a_{21}}{a_{11}-r_{11}}
    ,\quad
    t_{11} = 1 - \frac{a_{11}}{r_{11}}
    \;.
\end{equation}
Using $Q_{mp}^g$, on the other hand, results in a $Q_{mm}$ factor that is a Givens matrix, and cannot be represented in a compact WY form. Using eqs.~\eqref{eq:dQmp} and~\eqref{eq:Zpn} allows in both cases the calculation of $\dd Q_{mp}$, but the result is only correct when using $Q_{mp}^h$.

\end{document}